\documentclass{amsart}
\usepackage{amssymb}

\usepackage{mathtools}
\usepackage{mathrsfs}

\usepackage{xcolor}

\newtheorem{theorem}{Theorem}[section]
\newtheorem{lemma}[theorem]{Lemma}
\newtheorem{proposition}[theorem]{Proposition}

\theoremstyle{definition}

\theoremstyle{remark}

\makeatletter
\newcommand{\cp}{\mathop{\operator@font cp}}
\newcommand{\range}{\mathop{\operator@font range}}
\newcommand{\rank}{\mathop{\operator@font rank}}
\newcommand{\dom}{\mathop{\operator@font dom}}
\newcommand{\Real}{\mathop{\operator@font Re}}
\newcommand{\cont}{\mathop{\operator@font Cont_w}}
\newcommand{\alg}{\mathop{\operator@font Alg\mathcal{N}}}
\newcommand{\tr}{\mathop{\operator@font tr}}
\newcommand{\sgn}{\mathop{\operator@font sgn}}
\newcommand{\kn}{\mathop{\operator@font \mathcal{K}(\mathcal{N})}}
\newcommand{\scr}{\mathop{\operator@font C_0(X)\times_{\phi}\mathbb{Z}_+}}
\newcommand{\supp}{\mathop{\operator@font supp}}
\newcommand{\rad}{\mathop{\operator@font Rad}}
\newcommand{\hrad}{\mathop{\operator@font HRad}}
\newcommand{\D}{\mathop{\operator@font D}}
\newcommand{\hc}{\mathop{\operator@font \mathcal{R}_{hc}(\mathcal{A})}}
\@namedef{subjclassname@2020}{%
  \textup{2020} Mathematics Subject Classification}
\makeatother

\begin{document}

\title{Topological radicals of semicrossed products}
\author{G. Andreolas}
\author{M. Anoussis}
\author{C. Magiatis}
%\date{\today}
\address{Department of Mathematics, University of the Aegean, 832\,00
Karlovassi, Samos, Greece}
\email{gandreolas@aegean.gr}
\address{Department of Mathematics, University of the Aegean, 832\,00
Karlovassi, Samos, Greece}
\email{mano@aegean.gr}
\address{Department of Mathematics, University of the Aegean, 832\,00
Karlovassi, Samos, Greece}
\email{chmagiatis@aegean.gr}


\begin{abstract}
 We characterize the hypocompact radical of a semicrossed product in terms of properties of the dynamical system.
We show that an element $A$ of a semicrossed product is in the hypocompact radical if and only if the Fourier 
coefficients of 
$A$ vanish on the closure of the recurrent points and the $0$-Fourier coefficient vanishes also on the largest 
perfect subset of $X$. 
\end{abstract}

\maketitle 

%%%%%%%%%%%%%%%%%%%%%%%%%%%%%%%%%%%%%%%%%%%%%%%%%%%%%%%%%%%%%%%%%%%%%%%%%%%%%%%%%%%%%%%%%%%%%%%%%%%%%%%%%%%%%%%%%
\section{Introduction and Preliminaries}

 Let $\mathcal B$ be a Banach algebra.
 An element $a$ of  $\mathcal{B}$ is said to be compact  if the map $M_{a,a}:\mathcal{B}\rightarrow\mathcal{B}$, $x\mapsto axa$ is 
compact. 
Following Shulman and Turovskii \cite[3.2]{shutu4} we will call a  Banach algebra 
$\mathcal{B}$  \textit{hypocompact}  if any nonzero quotient $\mathcal{B}/\mathcal{J}$ by a closed ideal 
$\mathcal{J}$ contains a nonzero
compact element. We will say that an  ideal $\mathcal J$ of a Banach algebra $\mathcal B$ is  hypocompact if it is  
hypocompact  as an
algebra. Shulman and Turovskii have proved that any Banach  algebra $\mathcal{B}$ has a largest  hypocompact 
ideal \cite[Corollary 3.10]{shutu4}. This ideal is closed  and
 is called the hypocompact  radical of $\mathcal{B}$. We will  denote it  by $\mathcal{B}_{hc} $.

The hypocompact radical of Banach algebras was studied within the framework of the theory of topological radicals
\cite{shutu4, shutu5}.
This theory  originated 
with Dixon \cite{dix} and was further developed by  
Shulman and Turovskii 
 in a series of papers
\cite{shutu0,  shutu1, shutu2, shutu4, shutu5} and by  Kissin, Shulman and Turovskii \cite{shutu3}.
 The theory of topological radicals has applications  to  various problems of Operator Theory and 
Banach algebras.

 It follows from \cite[Lemma 8.2]{bre}, that the hypocompact radical contains the ideal generated by the compact elements.
If  $\mathcal X$ is a Banach space, we shall denote by $\mathcal{B}(\mathcal X)$ the Banach algebra of all 
bounded linear operators on $\mathcal X$ and by
$\mathcal{K}(\mathcal X)$ the Banach subalgebra of all compact operators on $\mathcal X$. 
Vala has shown in \cite{vala64} that  an element $a \in \mathcal{B}(\mathcal X)$ is a compact element if 
and only if $a \in \mathcal{K}(\mathcal X)$.   It follows that if 
   $\mathcal H$ is a separable Hilbert space,  the hypocompact radical of $\mathcal{B}(\mathcal H)$ is $\mathcal{K}(\mathcal H)$. 
 Indeed, the ideal  $\mathcal{K}(\mathcal H)$ is the only proper ideal of $\mathcal{B}(\mathcal H)$ while the Calkin algebra 
$\mathcal{B}(\mathcal H)/\mathcal{K}(\mathcal H)$ does not have any non-zero compact element \cite[section 5]{fs}.
 
  Shulman and Turovskii observe   in \cite[p. 298]{shutu4} that  there exist Banach spaces $\mathcal X$, such that the hypocompact radical
 $\mathcal{B}(\mathcal X)_{hc}$ of $\mathcal{B}(\mathcal X)$ contains all the weakly compact
 operators and  contains strictly 
 the ideal of compact operators $\mathcal{K}(\mathcal X)$.  

  Argyros and Haydon constructed in \cite{ah}   a Banach space $\mathcal X$ 
 such that every operator in $\mathcal{B}(\mathcal X)$ is  a scalar 
multiple of the identity plus a compact operator. It follows that $\mathcal{B}(\mathcal X)/\mathcal{
K}(\mathcal X)$ is finite-dimensional and hence the hypocompact radical of $\mathcal{B}(\mathcal X)$ coincides with 
$\mathcal{B}(\mathcal X)$. 

A nest $\mathcal{N}$ on a Hilbert space $\mathcal H$ is a totally
ordered family of closed subspaces of  $\mathcal H$ containing $\{0\}$ and $\mathcal H$, which is closed under intersection 
and
closed span. If  $\mathcal{N}$ is a nest on a Hilbert space $\mathcal H$,  the nest algebra associated to $\mathcal N$ 
is the (non selfadjoint) algebra of 
all
operators $T\in\mathcal{B}(\mathcal H)$ which leave each member of $\mathcal N$ invariant.
The hypocompact radical of a nest algebra was characterized in \cite{anan}.

We recall the construction of the semicrossed product we will consider in this work.
Let $X$  be a locally compact metrizable space and $\phi:X\rightarrow X$ a homeomorphism. The pair $(X, \phi)$ 
is called a dynamical system.
 An action 
of $\mathbb{Z}_+$ on $C_0(X)$ by isometric $*$-automorphisms $\alpha_n$, $n\in\mathbb{Z}_+$ is obtained by defining 
$\alpha_n(f)=f\circ\phi^n$.
We write the elements of the Banach space $\ell^1(\mathbb 
Z_+,C_0(X))$ as formal series $A=\sum_{n\in\mathbb Z_+}U^nf_n$ with the norm given by 
$\|A\|_1=\sum_{n\in\mathbb Z_+}\|f_n\|_{C_0(X)}$.  The multiplication on $\ell^1(\mathbb Z_+,C_0(X))$ is 
defined by setting
$$U^nfU^mg=U^{n+m}(\alpha^m(f)g)$$
and extending by linearity and continuity. With this multiplication, $\ell^1(\mathbb Z_+,C_0(X))$ is a Banach algebra.

The Banach algebra $\ell^1(\mathbb Z_+,C_0(X))$ can be faithfully represented as a (concrete) operator algebra on a 
Hilbert space. This is achieved by assuming a faithful action of $C_0(X)$ on a Hilbert space 
$\mathcal{H}_0$. Then, we can define a faithful contractive representation $\pi$ of $\ell_1(\mathbb Z_+,C_0(X))$ on the Hilbert 
space $\mathcal H=\mathcal{H}_0\otimes \ell^2(\mathbb Z_+)$ by defining $\pi(U^nf)$ as
$$\pi(U^nf)(\xi\otimes e_k)=\alpha^k(f)\xi\otimes e_{k+n}.$$
The \emph{semicrossed product} $C_0(X)\times_{\phi}\mathbb Z_+$ is the closure of the image of 
$\ell^1(\mathbb Z_+,C_0(X))$ in $\mathcal{B(H)}$ in the representation just defined, 
where $\mathcal{B(H)}$ is the algebra of bounded linear operators on $\mathcal{H}$. 
Note that the semicrossed product is in
fact independent of the faithful action of $C_0(X)$ on $\mathcal{H}_0$ (up to isometric
isomorphism) \cite{dkm}.
We will denote the semicrossed product $C_0(X)\times_{\phi}\mathbb Z_+$ by $\mathcal{A}$ and an
element $\pi(U^nf)$ of $\mathcal{A}$ by $U^nf$ to simplify the notation. 
We refer to \cite{pet, dkm, dfk}, for more information about the semicrossed product.

For $A=\sum_{n\in\mathbb Z_+}U^nf_n\in \ell^1(\mathbb Z_+,C_0(X))$, we call $f_n\equiv 
E_n(A)$ the $n$th \emph{Fourier coefficient} of $A$. The maps $E_n:\ell^1(\mathbb Z_+,C_0(X))\rightarrow C_0(X)$ are contractive in the (operator) norm of $\mathcal A$, and therefore they extend to
contractions $E_n:\mathcal A \rightarrow C_0 (X)$. An element $A$ of the semicrossed product $\mathcal A$  is $0$ if and only if 
$E_n(A)=0$
for all $n \in \mathbb Z+$ and thus $A$ is completely determined by its Fourier coefficients. We will denote $A$ by the formal series 
 $A=\sum_{n\in\mathbb Z_+}U^nf_n$, where $f_n=E_n(A)$. Note however that the series
$\sum_{n\in\mathbb Z_+}U^nf_n$ does not in general converge to $A$ \cite[II.9, IV.2 Remark]{pet}.

In this paper we characterize the hypocompact radical of a semicrossed product in terms of properties of the dynamical system.
We show that an element $A$ of a semicrossed product is in the hypocompact radical if and only if the Fourier coefficients of 
$A$ vanish on the closure of the recurrent points and the $0$-Fourier coefficient vanishes also on the largest 
perfect subset of $X$.

\section{The    hypocompact radical }

To obtain the characterization of the hypocompact radical of  a semicrossed product 
we recall the following properties of the hypocompact radical of a Banach algebra proved by  
Shulman and  Turovskii in \cite{shutu4}.

\begin{theorem}\label{st}
 Let $\mathcal B$ be a Banach algebra and $\mathcal I$ a closed ideal of $\mathcal B$.
\begin{enumerate}
 \item 
 If $\mathcal B$ is hypocompact, then $\mathcal I$ and $\mathcal B/\mathcal I$ are hypocompact \cite[Corollary 3.9]{shutu4}.
 \item 
 If $\mathcal I$ and $\mathcal B/\mathcal I$ are hypocompact, then $\mathcal B$ is hypocompact \cite[Corollary 3.9]{shutu4}.
 \item
 Let $p: \mathcal B \rightarrow \mathcal B/\mathcal I$ be the  quotient map.
  Then $p(\mathcal B_{hc})\subseteq (\mathcal B/\mathcal I)_{hc}$\cite[Corollary 3.13]{shutu4}.
 \end{enumerate}
\end{theorem}

Let  X be a locally compact metrizable space. We shall use the 
 characterization of the   hypocompact radical of $C_0(X)$ which may  be obtained  using  \cite[Corollary 8.19 \& Theorem 8.22]{shutu5}.
We provide a   proof for completeness. 

A point $x\in X$ is called \emph{accumulation point} of $X$, if $x\in\overline{X\setminus\{x\}}$.
The set  of the accumulation points of $X$ is denoted $X_a$.
If $x\in X\setminus X_a$, then the point $x$ is called an \emph{isolated point}.
A subset $Y$ of a topological space is said to be \emph{dense in itself},
if it contains no isolated points.
If $Y$ is closed and dense in itself, it is said to be a \emph{perfect set}.
The set $Y$ is said to be a \emph{scattered set}, if it does not contain dense in themselves subsets.

It is well known that every space is   the disjoint union of a
perfect and a scattered one, and this decomposition is unique \cite[Theorem 3, p 79]{kura}.
If $X$ is a locally compact metrizable space, we write $X=X_p\cup X_s$ where $X_p$ is  the perfect set and $X_s$ is  the scattered set.

\begin{theorem}\label{hccox} 
If X is a locally compact metrizable space, then
\begin{equation*}
C_0(X)_{hc}=\left\{f\in C_0(X):f(X_p)=\{0\}\right\}.
\end{equation*}
\end{theorem}
\begin{proof}

Let $\mathcal I$ be the ideal  $\left\{f\in C_0(X):f(X_p)=\{0\}\right\}$ of $C_0(X)$. The ideal 
$\mathcal I$ is isomorphic to $C_0(X_s)$. We show that  every non-zero quotient of $\mathcal I$ 
by a closed ideal has a non-zero compact element.
Let  $\mathcal J$ be   a closed ideal of $\mathcal I$ and $S$  a 
closed subset  of $X_s$ such that $\mathcal J=\{f \in C_0(X_s): f(S)=\{0\}\}$.
 The quotient algebra 
$\mathcal I/\mathcal J$ is isomorphic to 
$C_0(S)$. Hence it suffices to   prove that the algebra $C_0(S)$ has a non-zero compact element. Since  
the set $S$ is contained in $X_s$ it is scattered, and  it contains  an isolated point $y$. Let $\chi_{\{y\}}$ be the 
characteristic function of $\{y\}$.
Then, the operator $M_{\chi_{\{y\}}, \chi_{\{y\}}}: C_0(S)  \rightarrow C_0(S)$ is a rank-one operator and hence,  
$\chi_{\{y\}}$ is a compact   element 
of the algebra $C_0(S)$.  It follows that 
$\mathcal I\subseteq C_0(X)_{hc}$.

We show now  that $\mathcal I= C_0(X)_{hc}$. 
Assuming   that $\mathcal I \neq C_0(X)_{hc}$    we will prove that the quotient algebra $C_0(X)_{hc}/\mathcal I$ 
contains no non-zero compact elements. This implies that $\mathcal I= C_0(X)_{hc}$ by Theorem \ref{st}.
Let $f\in C_0(X)_{hc}\setminus \mathcal I$. There exists $x_p\in X_p$, such that $f(x_p)\neq 0$ and  an 
open neighborhood $U_p$ of $x_p$, such that
\begin{eqnarray*}
|f(x)|>\frac{|f(x_p)|}{2}, & \forall x\in U_p.
\end{eqnarray*}

Consider  a sequence of points $\{x_i\}_{i\in\mathbb N}\subseteq U_p\cap X_p$ and a sequence 
of open subsets $\{V_i\}_{i\in\mathbb N}$ of $X$, such that $x_i\in V_i\subseteq U_p$ and $V_i\cap V_j=\emptyset$ 
for $i\neq j$. 

By Urysohn's lemma there exists  a sequence of 
functions $\{h_i\}_{i\in\mathbb N}$ such that $h_i(x_i)=1$ and $h_i(X\setminus V_i)=\{0\}$.
Let $q:C_0(X)_{hc}\rightarrow C_0(X)_{hc}/\mathcal I$ be  the quotient map. 
We estimate for $i\neq j$:
\begin{eqnarray*}
\|M_{q(f), q(f)}(q(h_i))-M_{q(f), q(f)}(q(h_j))\| & = & \inf_{g\in \mathcal I}\|f^2h_i-f^2h_j+g\|\\
& \ge & \inf_{g\in \mathcal I}|(f^2h_i-f^2h_j+g)(x_i)|\\
& = & |f^2(x_i)|>\frac{|f(x_p)|^2}{4}.
\end{eqnarray*}
Hence,  the sequence $\{M_{q(f), q(f)}(q(h_i))\}_{i\in\mathbb N}$
has no convergent subsequence, which implies that the element $q(f)$ is non compact.
\end{proof}

 Recall that  a set $Y\subseteq  X$ is called wandering if the sets 
$\phi^{-1}(Y), \phi^{-2}(Y), ...$ are pairwise disjoint. Since  $\phi$ is a homeomorphism, this condition is equivalent to
the condition that $\phi^m(Y)\cap \phi^n(Y)=\emptyset$, for all $m, n \in \mathbb Z_+, m\neq n$. 
A point $x \in X$ is called wandering if it possesses an open wandering neighborhood. Otherwise it is called non wandering. 
We will denote by $X_w$ the set of wandering points of $X$. It is clear that $X_w$  is 
the the union of all open wandering subsets of $X$.

Let $X_{1}$ be the set of non wandering points of $X$ and set $\phi_1=\phi|_{X_1}$ the restriction of $\phi$ to $X_1$. 
We thus obtain a dynamical system
$(X_1, \phi_1)$.
Define by transfinite recursion 
 a family $(X_{\gamma}, \phi_{\gamma})$ of dynamical systems. If $(X_\gamma, \phi_\gamma)$  is defined, then set $X_{\gamma+1}$ 
  the set of non wandering points of the dynamical system 
$(X_\gamma, \phi_\gamma)$ and $\phi_{\gamma+1}=\phi|{X_{\gamma+1}}$. If $\gamma$ is a limit ordinal and the systems 
$(X_\beta, \phi_\beta)$ have been defined for all $\beta<\gamma$, 
set 
$X_\gamma=\cap_{\beta<\gamma} X_\beta$ and $\phi_{\gamma}=\phi|_{X_{\gamma}}$  the restriction of $\phi$ to $X_\gamma$. 
This process must stop 
at some ordinal  $\gamma_0$, since the cardinality of
the family cannot exceed the cardinality of the power set of $X$. The following is   \cite[Lemma 13]{dkm}.

\begin{proposition}\label{cen}
The set $X_{\gamma_0}$ is the closure of the set of recurrent points $X_r$ of the system $(X, \phi)$.
\end{proposition}

If $\gamma$ is an ordinal $\gamma\leq \gamma_0$, we will denote by  $\mathcal I_\gamma$  the ideal 
 $$\left\{A \in\mathcal A:E_0(A)=0, E_n(A)(X_{\gamma})=\{0\}, \forall n \in \mathbb Z_+, n\geq1\right\}.$$

The proof of the following lemma is straightforward, and is omitted.
\begin{lemma}\label{lemma}
 If $\gamma$ is a limit ordinal, then $I_{\gamma}=\overline{\cup_{\beta<\gamma}I_\beta}$.
\end{lemma}

It is known that the ideal generated by the compact elements of $\mathcal{A}$ is contained in the hypocompact radical \cite{bre}.
We will need the following characterization of this ideal  which is proved in \cite{aam}.
\begin{theorem}\label{cmo31}
 The ideal generated by the compact elements of $\mathcal{A}$ is 
 the set 
 $$
 \{A\in\mathcal{A}\; |\; E_n(A)(X\setminus X_w)=\{0\},\; \forall n\in\mathbb{Z_+}\ \text{and } E_0(A)(X_a)=\{0\}\}.$$
\end{theorem}

The following is the main result of the paper.

\begin{theorem}\label{main}
 The hypocompact radical 
 $\mathcal A_{hc}$ of  $\mathcal A$  is equal to
 $$\mathcal I=\left\{A \in\mathcal A:E_0(A)(X_p)=0, E_n(A)(X_{\gamma_0})=\{0\}, \forall n \in \mathbb Z_+\right\}.$$
 
\end{theorem}

\begin{proof}

\textbf{1st step}

We shall  prove that  $\mathcal I$ is contained in $\mathcal A_{hc}$.
We first prove that  $\mathcal I_{\gamma_0}$ is contained in $\mathcal A_{hc}$.
Assume the contrary.

 It follows from  Theorem \ref{cmo31} that $\mathcal I_1$  is contained in the ideal generated by the compact elements. 
  The hypocompact radical contains the ideal generated by the compact elements \cite{bre}, and hence $\mathcal I_1$ is contained 
 in $\mathcal A_{hc}$.

 Let  $\beta$ be the least ordinal $\beta\leq \gamma_0$ such that $\mathcal I_\beta$  is not contained in $\mathcal A_{hc}$.
 We show that $\beta$ is  a successor. If not,  since 
 $\mathcal I_{\gamma}\subseteq \mathcal A_{hc}$  for all $\gamma <\beta$, we obtain from Lemma \ref{lemma} that 
 $\mathcal I_{\beta}=\overline{\cup_{\gamma<\beta}\mathcal I_\gamma}\subseteq \mathcal A_{hc}$, which is absurde. Hence, 
$\beta$ is a successor.
 
 We are going to prove  that $\mathcal I_{\beta}$ is a hypocompact algebra.
 Consider the algebra $\mathcal I_{\beta}/\mathcal I_{\beta-1}$. It suffices to show that 
 $\mathcal I_{\beta}/\mathcal I_{\beta-1}$ is hypocompact, since the class of
 hypocompact algebras is closed under extensions and the ideal $\mathcal I_{\beta-1}$ is hypocompact (Theorem \ref{st}).
  
 We show that the algebra $\mathcal I_{\beta}/\mathcal I_{\beta-1}$ is generated by the compact elements it contains  
and hence is a hypocompact algebra 
 by  \cite{bre}.
 
 Let $A \in \mathcal I_{\beta}$. It follows  from the condition defining $\mathcal I_{\beta}$, that $U^nE_n(A) \in 
\mathcal I_{\beta}$,
 for all $n\in \mathbb Z_+, n\geq1$.
Hence, it suffices to show that the image of $U^nE_n(A)$ under  the 
natural map $\pi: \mathcal I_{\beta}\rightarrow \mathcal I_{\beta}/\mathcal I_{\beta-1}$ is contained in the ideal  
generated by the compact elements of 
$\mathcal I_{\beta}/\mathcal I_{\beta-1}$. It suffices to see this for an element of $\mathcal I_{\beta}$ of the form  
$U^nf$ with $f$ compactly supported.
 It follows from  \cite[Lemma 14]{dkm}, that $f$ can be written as a finite sum $f=\sum f_i$ where  each 
 $f_i$ has compact support contained in an open set $V_i$ such that 
 $ V_i \cap X_{\beta-1}$ is  wandering  for the system
 $(X_{\beta-1}, \phi_{\beta-1})$ and    $U^nf_i \in \mathcal I_{\beta}$, for all $i$.
 
 Hence, it suffices to prove that $\pi(U^nf)$ is a compact element, where 
 $f$  has compact support contained in an open set $V$, such that 
 $ V \cap X_{\beta-1}$ is  wandering for the system
 $(X_{\beta-1}, \phi_{\beta-1})$.

 We calculate:

 $$U^nf(\sum U^mg_m)U^nf=\sum U^{2n+m}f\circ \phi^{m+n}g_m\circ \phi^{n}f,$$
 
 for  $\sum U^mg_m \in \mathcal I_{\beta}$.
 
 Since $n\geq 1$, we have $n+m \geq 1$, for all $m \in \mathbb Z_+$, and consequently 
 $f\circ \phi^{m+n}f=0$ on $X_{\beta-1}$, for all $m \in \mathbb Z^+$ since  $V \cap X_{\beta-1}$ 
 is wandering. Hence, $U^nf( U^mg_m)U^nf= U^{2n+m}f\circ \phi^{m+n}g_m\circ \phi^{n}f \in \mathcal I_{\beta-1}$.

 Thus, $\pi(U^nf)$ is a compact element of $\mathcal I_{\beta}/\mathcal I_{\beta-1}$, and $\mathcal I_{\beta}$ is 
 a hypocompact ideal which is a contradiction. We conclude that $\mathcal I_{\gamma_0}$ is contained in 
 $\mathcal A_{hc}$.
Now, $\mathcal I/\mathcal I_{\gamma_0}$ is isomorphic to $\{f \in C_0(X): f(X_p\cup X_{\gamma_0})=\{0\}\}$ 
which is a hypocompact algebra by
Theorem \ref{hccox}.
It follows from Theorem \ref{st}  that $\mathcal I$ is a hypocompact ideal, and hence it is contained in  $\mathcal A_{hc}$.

\textbf {2nd step}

 We show now that 
 $\mathcal A_{hc}=\mathcal I$.
We will suppose that $\mathcal I\subsetneq \mathcal A_{hc}$ 
and we will prove that the quotient algebra $\mathcal A_{hc}/ \mathcal I$, contains no non-zero compact elements.
This  implies that $\mathcal A_{hc}=\mathcal I$ by Theorem \ref{st}.

Let $A\in\mathcal A_{hc}\setminus \mathcal I$ and set  $E_m(A)=f_m$, for all $m\in\mathbb{Z}_+$. 
Since the map $E_0$ is a continuous homomorphism from $\mathcal A$ onto $C_0(X)$, 
it follows
from Theorem \ref{st} that $E_{0}(\mathcal A_{hc})\subseteq C_0(X)_{hc}$ and hence 
by  Theorem \ref{hccox} we have $E_0(A)(X_p)=\{0\}$.

Since $A \notin \mathcal I$, it follows from Proposition \ref{cen} 
that there exists $m\in\mathbb Z_+$ such that $f_m(X_r)\neq \{0\}$. We set 
\begin{eqnarray*}
m_0=\min\{m\in\mathbb Z_+:  f_m(X_r)\neq \{0\}\},
\end{eqnarray*}
and we consider $x_0\in X_r$ such that $f_{m_0}(x_0)\neq 0$. There exists  an open neighborhood $U_0$ of $x_0$ such 
that 
\begin{eqnarray}\label{hcre0}
|f_{m_0}(x)|>\frac{|f_{m_0}(x_0)|}{2} & , & \forall x\in U_0.
\end{eqnarray}

Since  $x_0$ is a recurrent point,  there exist an open neighborhood 
$V_0$ of $x_0$ such that $\overline{V_0}\subseteq U_0$ and a  strictly increasing 
sequence $\{n_{i}\}_{i=1}^{\infty}\subseteq\mathbb N$ such that 
\begin{eqnarray}\label{hcre01}
\phi^{n_i}(x_0)\in V_0 & , & \forall i\in\mathbb N.
\end{eqnarray} 
Choosing, if necessary, a subsequence, we may assume that  $n_1> m_0$ and $n_{i+1}>3n_i$. 
By Urysohn's  lemma there is $u_0\in C_0(X)$ such that $u_0(x)=1$, for all 
$x\in \overline{V_0}$ and $u_0(X\setminus U_0)=\{0\}$.  We thus have
\begin{eqnarray}\label{hcre1}
u_0(x_0)=u_0\circ\phi^{n_i}(x_0)=1 & , & \forall i\in\mathbb N.
\end{eqnarray}

By \cite[Proposition 2.1]{muhly},
 we have that $U^{m_0}f_{m_0}\in\mathcal A_{hc}$, (see also \cite[p. 133]{dkm}). 
Hence, if we  consider the sequence $\{B_i\}_{i=1}^{\infty}$, where
\begin{eqnarray*}
B_i & = & (U^{n_{i+1}-n_i-m_0}u_0\circ\phi^{-m_0})(U^{m_0}f_{m_0})(U^{n_i-m_0}u_0\circ\phi^{-m_0})\\
& = & U^{n_{i+1}-m_0}u_0\circ\phi^{n_i-m_0}f_{m_0}\circ\phi^{n_i-m_0}u_0\circ\phi^{-m_0},
\end{eqnarray*}
it follows that
$\{B_i\}_{i=1}^{\infty}\subseteq\mathcal A_{hc}$.

Let $\pi:\mathcal A_{hc}\rightarrow \mathcal A_{hc}/\mathcal I$ be the  quotient  map. 
To prove that the element $\pi(A)$ is not a  compact element of 
$\mathcal A_{hc}/\mathcal I$, we will prove that the sequence $\{M_{\pi(A),\pi(A)}(\pi(B_i))\}_{i\in \mathbb N}$ 
has no Cauchy subsequence.

Let  $k,l\in\mathbb N$ with $k>l$. 
If $r< (n_{k+1}-m_0)$, the $r$th  Fourier coefficient of $B_k$ is $0$, and this also holds for
$M_{A,A}(B_k)$. It follows that 
\begin{eqnarray*}
E_{n_{l+1}+m_0}(M_{A,A}(B_k))=0,
\end{eqnarray*}
since $n_{l+1}+m_0<3n_{l+1}-m_0<n_{k+1}-m_0$.

Therefore,  it follows that 
\begin{eqnarray*}
\|M_{\pi(A),\pi(A)}(\pi(B_k-B_l))\| & = & \inf_{N\in \mathcal I}\|M_{A,A}(B_k-B_l)+N\|\\
& \ge & \inf_{N\in \mathcal I}\|E_{n_{l+1}+m_0}\left(M_{A,A}(B_k-B_l)+N\right)\|\\
& \ge & \inf_{N\in \mathcal I}|E_{n_{l+1}+m_0}\left(M_{A,A}(B_l)+N\right)(x_0)|\\
& = & |E_{n_{l+1}+m_0}(M_{A,A}(B_l))(x_0)|
\end{eqnarray*}
since $x_0 \in X_r$  and thus, for all  $N \in \mathcal I$, we have $E_{n_{l+1}+m_0}(N)(x_0)=0$.

We calculate   $|E_{n_{l+1}+m_0}(M_{A,A}(B_l))(x_0)|$.

We have 
\begin{eqnarray*}
&&|E_{n_{l+1}+m_0}(M_{A,A}(B_l))(x_0)|  = \\
&&\left|\sum_{n=0}^{2m_0}(f_{2m_0-n} \circ\phi^{n_{l+1}+n-m_0}u_0\circ\phi^{n_l+n-m_0}f_{m_0}\circ\phi^{n_l+n-m_0}u_0\circ\phi^{n-m_0}f_n)(x_0)\right|.
\end{eqnarray*}

For $n<m_0$ we have $f_n(x_0)=0$. Also, 
for $n>m_0$ and $n\leq 2m_0$ we have $f_{2m_0-n}\circ\phi^{n_{l+1}+n-m_0}(x_0)=0$, 
since $2m_0-n<m_0$ and $\phi^{n_{l+1}+n-m_0}(x_0) \in X_r$. 

Finally, 

\begin{eqnarray*}
&&|E_{n_{l+1}+m_0}(M_{A,A}(B_l))(x_0)|  = \\
&&|(f_{m_0}\circ\phi^{n_{l+1}}u_0\circ\phi^{n_l}f_{m_0}\circ\phi^{n_l}u_0f_{m_0})(x_0)|  \ge \frac{|f^3_{m_0}(x_0)|}{8}.
\end{eqnarray*}

It follows
that the sequence $\{M_{\pi(A),\pi(A)}(\pi(B_i))\}_{i\in \mathbb N}$ 
contains no Cauchy subsequence, and hence $\pi(A)$ is not a compact element of $\mathcal A_{hc}/ \mathcal I$.

\end{proof}

%%%%%%%%%%%%%%%%%%%%%%%%%%%%%%%%%%%%%%%%%%%%%%%%%%%%%%%%%%%%%%%%%%%%%%%%%%%%%%%%

\section{the scattered radical}

The following  are taken from  \cite[8.2]{shutu5}. 
A Banach  
  algebra  
 is called \textit{scattered} if the
spectrum of every element $a\in \mathcal{A}$ is finite or countable.  A   Banach algebra 
$\mathcal{A}$ has a largest scattered ideal denoted by
$\mathcal{R}_{s}(\mathcal{A})$. This ideal is closed and  
   is called the 
 scattered radical of $\mathcal{A}$ \cite[Theorem 8.10]{shutu5}.

Since all $C^*$-algebras are semisimple and their quotients are again $C^*$-algebras, 
it follows from  \cite[Theorem 8.22]{shutu5} that $C_0(X)_{hc}=C_0(X)_{s}$.

Donsig, Katavolos and Manousos proved in \cite{dkm} a  characterization of the Jacobson radical for more general semicrossed products. 
The next theorem follows from their result
\cite[Theorem 18]{dkm}.

\begin{theorem}\label{rad}

 The Jacobson 
radical of $\mathcal A$ coincides with the set of operators 

 $$
 \{A\in\mathcal{A}\; |\; E_n(A)(X_r)=\{0\},\; \forall n\in\mathbb{Z_+}\ \text{and } E_0(A)=0\}.$$
\end{theorem}

It follows from Theorem \ref{main}  and the above characterization, that the Jacobson radical of $\mathcal A$ is contained in $\mathcal A_{hc}$.
Hence, from  \cite[Theorem 8.15]{shutu5} we obtain the following.

\begin{theorem}\label{xaaa}

\begin{eqnarray*} 
\mathcal A_{\mathrm{hc}}=  \mathcal A_s.
\end{eqnarray*}
\end{theorem}

\vspace{1em}

\noindent \textbf{Acknowledgements.} This research is co-financed by Greece and the European Union (European Social Fund- ESF) through the Operational Programme ``Human Resources Development, Education and Lifelong Learning 2014-2020'' in the context of the project ``Compactness Properties and Topological Radicals of Non-Selfadjoint Operator Algebras'' (MIS 5048197).

%%%%%%%%%%%%%%%%%%%%%%%%%%%%%%%%%%%%%%%%%%%%%%%%%%%%%%%%%%%%%%%%%%%%%%%%%%%%%%%%%%%%%%%%%%%%%%
\bibliographystyle{amsplain}

\end{document}